# Optimal run time for an EPQ model with scrap, rework, setup time and machine breakdown (failure) under learning effect assumption

Masoud Fekri[1]

In manufacturing systems, setup time and production time play important parts which are reduced through human learning effect. Furthermore, defective products of the manufacturing lines are classified as imperfect (repairable) and scraps (non-repairable) products. To this end, purpose of this research is to determine the optimal run time for an Economic Production Quantity (EPQ) model with scrap, rework, setup time, and machine breakdown under learning effect assumption. Machine breakdown occurs after generating a specific quantity of defective items which follows by reworking process. In this study a percentage of defective products are considered as scrap and the rest of them are reworked, while times of production and setup are reduced by learning effect. The proposed problem is extended upon a production-inventory cost function. Next, convexity of the developed model is proved, and the optimal run time for the model is obtained. Finally, numerical example illustrates utility of the developed model.

## 1. Introduction

In today's highly competitive world markets, companies need to employ an economical inventory policy to satisfy customer demand. In this regard, many researchers have developed inventory/production models. The Economic Production Quantity (EPQ) model as an extension of the Economic Order Quantity (EOQ) model which dates back to (Harris 1913) has been developed by (Taft in 1918). The classical production quantity is one of the most popular and appealing issues of the inventory management upon which the quantity of a company is determined in an economic manner to minimize the total inventory costs by balancing the inventory holding cost and the average fixed ordering cost (Hsieh and Dye 2012). However a considerable amount of production-inventory models with more complicated and/or practical assumptions have been extensively studied today and during past decades, the EPQ model has restrictions and cannot be regarded as a universal inventory model (Chiu et al. 2010, Taleizadeh *et al*. 2014). Since many key assumptions of the classical EPQ model relax the factors which are inevitable in many manufacturing systems.

[1] School of Industrial Engineering, Iran University of Science and Technology, Tehran 16846-13114
corresponding author: masoud_fekri@ind.ius.ac.ir

One of the key assumptions of the basic EPQ model is that many items are produced with the best quality, while defective items are inevitably resulted due to some reasons such as imperfect raw materials, low skill levels of the workers, inadequate machine capability, and poor maintenance policies or when the manufacturing process gets out of control, among others (Pentico and Drake 2011, Pasandideh *et al*. 2013, Taleizadeh *et al*. 2014, Sivashankari and Panayappan, 2014). Therefore, numerous researches have been studied to enhance the classic finite production rate model by addressing the issue of defective items. In practical production environments, the defective/imperfect products can be reworked and repaired, leading to overall production cost reduction (Chiu *et al*. 2010). In addition to the random defective rate, the breakdown of production equipment has been introduced in the literature body as another critical reliability factor that can be very disruptive when happening (Chiu et al. 2007, Chiu et al. 2010, Pasandideh et al. 2013, Taleizadeh *et al*. 2014). Thus, some inevitable features of the EPQ model may be categorized into several classes such as EPQ model with imperfect items, random breakdown, reworking the defective products (completely or partially), scrap, backlogging (e.g. Chiu 2003, Chiu *et al*. 2007); EPQ model with random breakdowns, stochastic demand and deterioration rate (e.g. Widyadana and Wee 2012); random breakdowns with resuming or not resuming the interrupted lots after failure; optimal production run length; EPQ model with learning effect theory (e.g. Hou 2007, Darwish 2008, Jaber and El Saadany 2011).

In this study, we attempt to cope with the shortcomings of the traditional EPQ models. To this end, an inventory-production model is proposed regarding to imperfect product, reworking of defective items, scrap, setup time, and machine breakdown/failure under learning effect assumptions. It must be noted, setup time and production time of the model are affected by the learning effects and it can be reduced in the next stages of the production and reworks of defective products. As previous works, it is supposed a random breakdown occurs during the production length. In spite of that, learning effect theory is considered, which leads to reduce production time, and it could be the reason of the convexity of inventory plot in both stages of production and reworking defective products. In other word due to the learning effects, the setup time of reworking stage becomes lower than the setup time of production stage whenever a breakdown occurs, and the learning effect cause the concavity of the plot. Furthermore, both stages of the production and reworking are concave in inventory of the defectives plot. Moreover in order to develop operational measures of the proposed EPQ model, a cost structure is defined which captures some different operational cost factors including production cost, breakdown and setup cost, reworking cost, scrap cost, holding cost upon which run time is optimized.

Rest of the paper is structured as follows. Section 2 briefly reviews the previously conducted researches. Section 3 presents the proposed EPQ model. And, the proposed solution methodology is analysed thoroughly in Section 4. The paper follows in Section 5 on numerical experiment to validate the proposed model. Finally, some conclusions and future research directions are provided in the last section.

## 2. Literature review

One of the main research topics in production and inventory management is EPQ model which was firstly introduced under various assumptions to obtain optimal quantity of produced items. Since then, researchers have extended the classical EPQ in numerous different directions by relaxing one or more of its assumptions to overcome its shortcoming (Tai 2013, Chiu *et al*. 2007).

One basic assumption in traditional EPQ model is that the produced items are of perfect quality (Tai 2013, Chiu *et al*. 2007). However, defective items are inevitably resulted due to the deterioration of process or many other factors in reality (Pentico and Drake 2011). Hence, a considerable amount of study has been extended to enhance the classic EPQ model by considering the issue of imperfection items produced (e.g. Lee 1992, Rezaei and Davoodi 2008, Yassine *et al*. 2012, Wee *et al*. 2013). Rezaei and (Davoodi 2008) considered a supply chain with multiple products and multiple suppliers subject to capacity limitation so that received items from suppliers were not of perfect quality. (Yassine *et al*. 2012) presented disaggregating the shipments of imperfect quality items in a single production run in which all imperfect quality items are detected by the end of the production cycle and aggregating the shipments of imperfect items over multiple production runs. (Wee *et al*. 2013) suggested an EPQ model with imperfect quality items, shortage backordering and screening constraint, and applied the Renewal Reward Theorem (RRT) to formulate the exact expected total profit per unit time.

In practical production milieus, another important issue is how to handle imperfect quality items. One possible way is to rework and repaired leading to decrease significantly overall production-inventory costs (Tai 2013, Chiu *et al*. 2007). Some researches of rework process of products were addressed by (Salameh and Jaber 2000,

Jaber and Bonney 2003) and (Stewart *et al*. 2004). (Sarker *et al*. 2008) developed Economic batch quantity models in a multi-stage rework system to minimize the total cost considering two operational policies in a multi-stage production system. Sarkar and (Moon 2011) have extended the classical EPQ model with stochastic demand under the effect of inflation and backordering in an imperfect production system produces a single item type in which some products are defective and can be reworked at a cost. (Tai 2013) derive an approximated optimal cycle length and total operating cost per cycle for two EPQ models for deteriorating/imperfect quality items with imperfect quality items in which the inspection process for deteriorating items was assumed to be imperfect as well. (Pasandideh *et al*. 2010) developed a multi-product EPQ model with rework imperfect items of different product types and a warehouse space-limitation. And, the authors proposed a genetic algorithm to solve the developed non-linear integer-programming model.

Additionally, there are some recent works based on dynamical approaches where the underlying system is modelled using System Dynamics (Sterman 1994, Rafieisakhaei *et al.* 2016 a, b). This method models the systems using stock and flow and casual loop casual diagrams. These tools helps the model to capture the uncertainty in the system and gives better understanding of the complex dynamic problems which could arise in different disciplines such as economical problems( Taleizadeh *et al*. 2015d)

Furthermore in classical EPQ model, imperfect products remain as surplus or scrap at the end of each production cycle. Therefore, segregating inappropriate products could be found as one of manners to develop the basic one. Thus, (Hayek and Salameh 2001) derived an optimal operating policy for the finite production model of imperfect quality items with reworking so that all defective items are repairable and allowed backorders. And then, (Chiu 2003) extended (Hayek and Salameh 2001) model by addressing scrap items and the reworking of repairable items with a random defective rate. And, an EPQ inventory model with scrapped products and limited production capacity was introduced (Taleizadeh *et al*. 2010a), and an EPQ inventory model has been developed by addressing production quantity with defective items, service level constraints, random breakdown, production capacity limited, and partial backordering (Taleizadeh *et al*. 2010c). In the case of multi-product single-machine production system, (Taleizadeh *et al*. 2010f) developed an EPQ model with stochastic scrapped production rate, partial backordering and service level constraint. (Chiu *et al*. 2011) developed a problem to incorporate rework process and multiple shipments into an imperfect EMQ model with random defective rate. In their model, the reworking of all defective items takes place after completing the regular production process in each cycle so that a portion of reworked items fails during the reworking and becomes scrap. And, an EPQ inventory model with rework process (Taleizadeh *et al*. 2012g). Besides, (Taleizadeh *et al*. 2013b) presented an economic production quantity model with repair failure and limited capacity. In a subsequent article, (Taleizadeh *et al*. 2014c) extended an EPQ inventory model with scrap, reworking of defective products, interruption in process, backordering in a multi products-single machine system. (Pasandideh *et al*. 2013) considered a bi-objective EPQ model with imperfect items, rework and scrap subject to order limitation to minimize total inventory cost as well as to minimize the required warehouse space for raw material and perfect products. Moreover, (Sivashankari and Panayappan 2014) studied the optimal lot-sizing decision for a single product manufactured in a single stage manufacturing system with planned backorders and that generates imperfect quality products. In this model, all defective items are reworked and where scrap is produced, detected and discarded during the reworking.

Besides the dynamical approaches, mathematical modelling also helps to fundamentally understand the system (Azadeh *et al*.2016, Barazandeh *et al*.2017). This approach is mostly applicable for systems with high level of certainty and helps to diagnose and detect faults in the system in order to improve the overall performance and efficiency (Ben-Daya *et al*.2012, Bastani *et al*.2018).

Another crucial reliability factor is the failure of the production equipment that can be very disruptive when happening. To this end, (Liu and Cao 1999) analysed a production-inventory model with Poisson demand and random machine breakdowns. (Chung 1997) investigated the upper and lower bounds of the optimal lot sizes in an EPQ model with machine breakdown. (Makis and Fung 1998) introduced an optimal production/inspection policy to study the effect of machine failures on the optimal lot size and on the optimal number of inspections in a production cycle on an EMQ model. (Kim and Hong 1999) studied the optimal production run length in deteriorating production processes with an arbitrary probability distribution. (Chung 2003) developed approximations to production lot sizing for EPQ model with random machine breakdowns. (Lin and Kroll 2006) discussed the effect of an imperfection in the production process (quality issues) and equipment (machine breakdowns) on the EMQ, in which all defective units can be reworked that the repair activity is completed before the end of a production cycle. (Chiu *et al*. 2007) presented an optimal production run time of EPQ model with scrap, reworking of defective items, and stochastic machine breakdowns under no-resumption (NR) inventory control policy. Then after, (Chiu and Ting 2010) re-examined a theorem on conditional convexity of the integrated total cost function of (Chiu *et al*. 2007) paper and presented a direct proof to the convexity of the total cost function that can be applied in place of it to enhance quality of their optimization process. In addition, (Chiu *et al*. 2010) presented an optimal production run time for an imperfect finite production rate model with scrap, rework, and stochastic machine breakdown under the abort/resume (AR) policy despite (Chiu *et al*. 2007) study which has been developed based on NR policy. In other study, (Ting and Chung 2013) considered the economic production quantity model with scrap, rework and stochastic machine breakdown in which the rigorous methods have been adopted to demonstrate that the expected total cost per unit time is convex on all positive numbers to improve the conditional convexity in Theorem 1 of (Chiu *et al*. 2010). Moreover, they derived bounds for the optimal production run time to remove the logical shortcomings of mathematics presented in proof of Theorem 2 of (Chiu *et al*. 2010). (Chiu and Chang 2014) later made an extension of (Chiu *et al*.2007) model to incorporate failure in rework factor to enhance the practical usages of their research outcome. In the study, the optimal replenishment run time for the amended model has been derived, and the effect of the failure in rework on the optimal operating policy has been investigated through the exact model. In addition, an EPQ model by considering scrap, rework and machine breakdown maintenance has been presented by (Liao *et al*. 2009), and another EPQ model that considered production shift from an in-control state to an out-of-control state in which the production facility may break down at any random point during the production uptime period by (Chakraborty *et al*. 2009). And, (Chung *et al*. 2011) derived an EPQ model for deteriorating items under machine unavailability and shortage assumption.

As discussed earlier in the literature body until now no research is done in the inventory-production systems with machine breakdown/failure under learning effect assumption. Therefore, this paper develops an EPQ inventory model with defective items, rework, scrap, setup time, and machine breakdown/failure under learning effect assumption.

The earliest learning effect representation is a geometric progression that expresses the decreasing cost/ time required to accomplish any repetitive operation. The theory states that as the total quantity of units produced doubles, the cost per unit declines by some constant percentage (Wright 1936). Learning curve is a power function formulation that is represented by

$$q_x = q_1 x^{-b} \qquad (1)$$

where $q_x$ is the time to produce the $x^{th}$ unit, $q_1$ the time to produce the first unit, $x$ the production count, and $b$ the learning effect exponent. In practice, the $b$ parameter value is often replaced by another index number that has had more intuitive appeal. This index is referred to as the "fix learning rate in production ($\alpha_1$), fixed learning rate in produced items is assumed to be defective ($\alpha_2$) and fixed learning rate in setup ($\alpha_s$)" which occurs each time the production output is doubled, where $\alpha_1, \alpha_2, \alpha_s = \frac{q_{2x}}{q_x} = \frac{q_1(2x)^{-b}}{q_1(x)^{-b}} = 2^{-b}$ and $b = -\log("learning rate")/\log(2)$. The time to produce $x$ units by integrating Eq. (1) over the proper limits is given as:

$$t_x = \sum_{n=1}^{x} q_1 n^{-b} \cong \int_0^x q_1 n^{-b} dn \qquad (2)$$

## 3. Problem definition

In the present study, a single-stage manufacturing system is proposed that produces single product with imperfect quality and machine breakdown under no resumption (NR) inventory control policy. The production line of the proposed manufacturing system has made up of five operations including the initial machine setup process, production process leading to generating perfect, and imperfect items which are reworked, second setup process after machine failure occur, reworking process, and finally customer demand estimation/presentation. Furthermore, learning effects have been proposed lead to reduction of setup time, production time, and reworking time as well. It is noted, machine breakdown take place when $x$ portion of produced items is defective. A portion of the defective items are repairable and the rest of them are scrap. Moreover, it is assumed that there is no interruption during the rework process. During the regular production uptime, $x$ portion of produced items is assumed to be defective and is generated at a production rate $\lambda$. Among the defective items, a $\theta$ portion is supposed to be scrap and the other portion can be reworked and repaired. In each production cycle time, after machine setup and regular production process with learning effects, all repairable defective items are reworked at a rate $P_1$ when the machine is repaired under learning effects assumption as well. The production rate $P$ is constant and is much larger than the demand rate $D$ which is deterministic and constant. The production rate of defective items $\lambda$ could be expressed as the production rate times the defective rate, $\lambda = Px$. Furthermore, production rate of perfect quality items must always be greater than or equal to sum of the demand rate and the production rate of defective items ($P - D - \lambda \geq 0$ or $P - \frac{D}{P} - x \geq 0$). In addition, it is assumed there are no shortages in this model.

The following notations are used thorough the paper to develop the proposed inventory model:

*Parameters*

$P$ : production rate,

$P_1$ : production rate of defective production is constant and larger than the demand rate D,

$D$ : demand rate,

$\lambda$ : defective rate during the regular production,

$x$ : portion of produced items is assumed to be defective,

$\theta$ : portion of reworked items is assumed to be scrap,

$h$ : holding cost for a produced unit,

$h_1$ : holding cost for each reworked item,

$c$ : unit production cost,

$c_r$ : repairing cost of produced items is assumed to be defective,

$c_s$ : disposal cost of reworked items is assumed to be scrap,

$\alpha_1$ : fixed learning rate in production, $0 < \alpha_1 \leq 1$

$\alpha_2$ : fixed learning rate in produced items is assumed to be defective, $0 < \alpha_2 \leq 1$

$\alpha_s$ : fixed learning rate in setup, $0 < \alpha_s \leq 1$

$b_{p1}$ : production learning coefficient,

$$b_{p1} = -\log \alpha_1 / \log 2 \qquad (3)$$

$b_{p2}$ : learning coefficient associated with produced items is assumed to be defective,

$$b_{p2} = -\log \alpha_2 / \log 2 \qquad (4)$$

$b_s$ : learning coefficient associated with setup,

$$b_s = -\log \alpha_s / \log 2 \qquad (5)$$

$M$ : repair cost of machine breakdown,

$K$ : setup cost of machine,

*Auxiliary variables*

$H_1$ : the level of on-hand inventory when machine breakdown occurs,

$H_2$ : the maximum level of on-hand inventory when machine is repaired and the reworking of defective items is completed,

$H_3$ : the level of on-hand inventory when machine is repaired and restored,

$tp_1'$ : production time before a random breakdown occurs without learning effect,

$tp_2'$ : time needed for reworking of defective items when machine breakdown takes place without learning effect,

$tp_2$ : time needed for reworking of defective items when machine breakdown takes place with learning effect,

$ts_1$ : time required to set up production,

$ts_2$ : time required for repairing the machine,

$td$ : time needed for consuming all available perfect quality items,

$T$ : cycle length,

PC: the total production cost,

BC: the total breakdown and setup cost,

RC: the total reworking cost,

SC: the total scrap cost,

HC: the total holding cost,

*Decision variables*

$tp_1$ : production time before a random breakdown occurs with learning effect,

$TC(tp_1)$ : the total inventory cost per cycle,

The objective of this research is to obtain the optimal solution for the above defined decision variables that minimizes total cost including production cost, breakdown and setup cost, reworking cost, scrap cost, holding cost:

$$TC(tp_1) = PC + BC + RC + SC + HC \tag{6}$$

## 4. Modelling and analysis

According to Figure 1, working cycle time consists of setup time, production time, and setup time for the reworking stage of defective products, reworking time of defective products and the consuming time of the non-defective products.

$$T = [ts_1 + ts_2] + [tp_1 + tp_2] + td \rightarrow td = T - [ts_1 + ts_2] - [tp_1 + tp_2] \tag{7}$$

The working cycle time is calculated by the following formula:

$$T = \frac{Ptp_1(1-\theta x)}{D} \tag{8}$$

Relation between production time before a random breakdown occurs with learning effect and time needed for reworking of defective items when machine breakdown takes place with learning effect is like Equation (9).

$$tp_2 = \frac{\lambda tp_1(1-\theta)}{P_1} \tag{9}$$

So the following formula can be obtained:

$$td = \frac{Ptp_1(1-\theta x)}{D} - [ts_1 + ts_1(2)^{-b_s}] - [tp_1 + \frac{\lambda tp_1(1-\theta)}{P_1}] \tag{10}$$

From Figure 1, one can obtain the level of on-hand inventory when machine breakdown occurs $H_1$, the level of on-hand inventory when machine is repaired and restored $H_3$ and the maximum level of on-hand inventory when machine is repaired and the reworking of defective items is completed $H_2$ as follows:

$$tp_1^{'} = \frac{H_1}{P - D - \lambda} \qquad (11)$$

$$tp_2^{'} = \frac{H_2 - H_3}{P_1 - D} \qquad (12)$$

we employ the relation between production time before a random breakdown occurs without learning effect and production time before a random breakdown occurs with learning effect likes Equation (13) and also Relation between time needed for reworking of defective items when machine breakdown takes place without learning effect and time needed for reworking of defective items when machine breakdown takes place with learning effect likes Equations (14):

$$tp_1^{'} - tp_1 = tp_1^{'} - tp_1^{'}(P - D - \lambda)(\int_0^{H_1} q^{-b_{p1}} dq) \rightarrow$$

$$tp_1 = tp_1^{'}(\int_0^{H_1} q^{-b_{p1}} dq) = \frac{H_1}{P - D - \lambda}(P - D - \lambda)\frac{1}{1 - b_{p1}}(H_1)^{1 - b_{p1}} \rightarrow tp_1 = \frac{(H_1)^{2 - b_{p1}}}{1 - b_{p1}} \qquad (13)$$

$$tp_2^{'} - tp_2 = tp_2^{'} - tp_2^{'}(P_1 - D)(\int_0^{H_1} q^{-b_{p1}} dq) \rightarrow$$

$$_2 = tp_2^{'}(P_1 - D)(\int_0^{H_1} q^{-b_{p1}} dq) = \frac{H_2 - H_3}{P_1 - D}(P_1 - D)\frac{1}{1 - b_{p2}}(H_2 - H_3)^{1 - b2} \rightarrow tp_2 = \frac{(H_2 - H_3)^{2 - b2}}{1 - b_{p2}}$$

$$(14)$$

Time required to setup production and time required for repairing the machine is as follows:

$$ts_2 = ts_1(2)^{-bs} \qquad (15)$$

Please insert Figure 1 here.

Please insert Figure 2 here.

As depicted in Figures 1 and 2, total production-inventory cost is:

$$TC(tp_1) = [c.P.tp_1] + [K + M] + [c_r.P.tp_1.x.(1-\theta)] + [c_s.P.tp_1.x.\theta] +$$

$$h\{[\frac{1}{P-D-\lambda}\int_0^{tp_1} t^{-\frac{1}{b_{p1}}}dt] + [\frac{1}{\lambda}\int_0^{tp_1} t^{\frac{1}{1-b_{p1}}}dt] + [H_3.ts_1(2)^{-b_s} + \frac{1}{D}\int_0^{D.ts_1(2)^{-b_s}} t^{-\frac{1}{1-b_s}}dt]$$

$$+[ts_1(2)^{-b_s}.\lambda.tp_1] + [H_3.tp_2 + \frac{1}{P_1-\lambda}\int_0^{tp_2} t^{-\frac{1}{b_{p2}}}dt] + [\frac{1}{2}H_2.td] + [h_1\frac{1}{P_1}\int_0^{tp_2} t^{-\frac{1}{1-b_{p2}}}dt]$$

$$\rightarrow TC(tp_1) = [c.P.tp_1] + [K + M] + [c_r.P.tp_1.x.(1-\theta)] + [c_s.P.tp_1.x.\theta] +$$

$$h\{[\frac{1}{P-D-\lambda}\frac{b_{p1}}{-1+b_{p1}}tp_1^{\frac{-1+b_{p1}}{b_{p1}}}] + [\frac{1}{\lambda}\frac{1-b_{p1}}{2-b_{p1}}tp_1^{\frac{2-b_{p1}}{1-b_{p1}}}] + [(^{2-b_{p1}}\sqrt{tp_1(1-b_{p1})} - D.ts_1(2)^{-b_s})ts_1(2)^{-b_s} +$$

$$[\frac{1}{D}\frac{b_s-1}{b_s}(D.ts_1(2)^{-b_s})^{\frac{b_s}{b_s-1}}] + [ts_1(2)^{-b_s}.\lambda.tp_1] + [(^{2-b_{p1}}\sqrt{tp_1(1-b_{p1})} - D.ts_1(2)^{-b_s})\frac{\lambda tp_1(1-\theta)}{P_1} +$$

$$\frac{1}{P_1-\lambda}\frac{b_{p2}}{1-b_{p2}}(\frac{\lambda tp_1(1-\theta)}{P_1})^{\frac{1-b_{p2}}{b_{p2}}}] + [\frac{1}{2}(^{2-b_{p2}}\sqrt{\frac{\lambda tp_1(1-\theta)}{P_1}}(1-b_{p2}) + ^{2-b_{p1}}\sqrt{tp_1(1-b_{p1})} - D.ts_1(2)^{-b_s})$$

$$(\frac{Ptp_1(1-\theta x)}{D} - (ts_1 + ts_1(2)^{-b_s}) - (tp_1 + \frac{\lambda tp_1(1-\theta)}{P_1}))]\} + h_1\frac{1}{P_1}\frac{b_{p2}-1}{b_{p2}}(\frac{\lambda tp_1(1-\theta)}{P_1})^{\frac{b_{p2}}{b_{p2}-1}}$$

(16)

Here we are primarily derived from the production-inventory cost function to find *TC (tp₁)* obtains a minimum:

$$\frac{dTC(tp_1)}{dtp_1} = cP + c_r Px(1-\theta) + c_s Px\theta + h\{\frac{1}{(P-D-\lambda)} tp_1^{-\frac{1}{b_{p1}}} + \frac{1}{\lambda} tp_1^{\frac{1}{1-b_{p1}}} +$$

$$\frac{(1-b_{p1})^{\frac{1}{2-b_{p1}}}}{2-b_{p1}} ts_1(2)^{-b_s} tp_1^{\frac{1}{2-b_{p1}}-1} + \frac{\lambda(1-\theta)}{P_1}(\frac{3-b_{p1}}{2-b_{p1}})(1-b_{p1})^{\frac{1}{2-b_{p1}}} tp_1^{\frac{3-b_{p1}}{2-b_{p1}}-1} - Dts_1(2)^{-b_s}\frac{\lambda(1-\theta)}{P_1} +$$

$$\frac{1}{P_1-\lambda}(\frac{b_{p2}}{1-b_{p2}})(\frac{\lambda(1-\theta)}{P_1})^{\frac{1-b_{p2}}{b_{p2}}} tp_1^{\frac{1-b_{p2}}{b_{p2}}-1} + \frac{1}{2}[\frac{P(1-\theta x)}{D} - 1 - \frac{\lambda(1-\theta)}{P_1}] \times$$

$$[(\frac{3-b_{p2}}{2-b_{p2}})(\frac{\lambda(1-\theta)}{P_1}(1-b_{p2}))^{\frac{1}{2-b_{p2}}} tp_1^{\frac{3-b_{p2}}{2-b_{p2}}-1} + (\frac{3-b_{p1}}{2-b_{p1}})(1-b_{p1})^{\frac{1}{2-b_{p1}}} tp_1^{\frac{3-b_{p1}}{2-b_{p1}}-1}] - \frac{1}{2}[ts_1 + ts_1(2)^{-b_s}] \times$$

$$[(\frac{1}{2-b_{p2}})(\frac{\lambda(1-\theta)}{P_1}(1-b_{p2}))^{\frac{1}{2-b_{p2}}} tp_1^{\frac{1}{2-b_{p2}}-1} + (\frac{1}{2-b_{p1}})(1-b_{p1})^{\frac{1}{2-b_{p1}}} tp_1^{\frac{1}{2-b_{p1}}-1}]\} + \frac{h_1}{P_1} tp_1^{-\frac{1}{1-b_{p1}}}$$

(17)

Also, because of the exponential and difficult and complicated integral relations in this function, it is difficult or impossible to acquire minimum of function. In this circumstance, we should use numerical values to the parameters to solve the problem. Also to demonstrate the convexity of the production-inventory cost function, second derivative of the function *TC (tp₁)* obtains and the result is as follow: *(0≤bp₁,bp₂≤1)*

Second derivation of the production-inventory cost function is composed of three words by negative value or four words by negative value *(if 0≤bp₂≤1/2 function is composed of three words by negative value and if 1/2≤bp₂≤1 function is composed of four words be negative value)* and other words with positive value (course with an assumption that will be listed below the that guarantees positive values for the word). The second derivative of the function must be positive to guarantee the convexity of the production-inventory cost function, otherwise we suppose impossible answer to this function. The minimum point of cost function should be find with respect to a particular *tp₁* so the second derivative of the function must be positive (If the Second derivation of function is negative, the answer is impossible). Therefore, numerical values that we should allot to problem parameters (here our targets is values of learning rate because value of other parameters in problem has been considered constant) should be chosen in a way that Second derivation of function be positive (answer should be possible).

Negative section in the second derivative of the function TC ($tp_1$):

1.
$$-h\frac{1}{(P-D-\lambda)}(\frac{1}{b_{p1}})tp_1^{-\frac{1}{b_{p1}}} \leq 0 \tag{18}$$

Because of this negative word,

$$(0 \leq b_{p1} \leq 1) \Rightarrow (\frac{1}{2-b_{p1}}-1) \leq 0 \tag{19}$$

entire term is negative:

2.
$$h\frac{(1-b_{p1})^{\frac{1}{2-b_{p1}}}}{2-b_{p1}}(\frac{1}{2-b_{p1}}-1)ts_1(2)^{-b_s} tp_1^{\frac{1}{2-b_{p1}}-2} \leq 0 \tag{20}$$

3.
$$\frac{h_1}{P_1}(-\frac{1}{1-b_{p1}})tp_1^{-\frac{1}{1-b_{p1}}-1} \leq 0 \tag{21}$$

And the fourth negative term is as follows:

$$(P_1 > \lambda) \Rightarrow \frac{1}{P_1 - \lambda} \geq 0$$
$$(0 \leq b_{p2} \leq 1) \Rightarrow \frac{b_{p2}}{1-b_{p2}} \geq 0 \tag{22}$$

4.
$$\text{if } (\frac{1-b_{p2}}{b_{p2}}-1 \leq 0 \rightarrow 1-b_{p2} \leq b_{p2} \rightarrow b_{p2} \geq \frac{1}{2}) \rightarrow$$
$$h(\frac{b_{p2}}{1-b_{p2}})(\frac{1-b_{p2}}{b_{p2}}-1)(\frac{\lambda(1-\theta)}{P_1})^{\frac{1-b_{p2}}{b_{p2}}} tp_1^{\frac{1-b_{p2}}{b_{p2}}-2} \leq 0 \tag{23}$$

On the other hand:

$$0 \leq b_{p2} \leq 1 \rightarrow \frac{1}{2} \leq b_{p2} \leq 1$$

Positive section in the second derivative of the function TC (tp1):

$$h\frac{1}{\lambda}(\frac{1}{1-b_{p1}})tp_1^{\frac{1}{1-b_{p1}}} \geq 0 \tag{24}$$

$$0 \leq b_{p1} \leq 1 \rightarrow \frac{3-b_{p1}}{2-b_{p1}} \geq 1 \rightarrow h\frac{\lambda(1-\theta)}{P_1}(\frac{3-b_{p1}}{2-b_{p1}})(\frac{3-b_{p1}}{2-b_{p1}}-1)(1-b_{p1})^{\frac{1}{2-b_{p1}}} tp_1^{\frac{3-b_{p1}}{2-b_{p1}}-2} \geq 0 \tag{25}$$

With this assumption that value of this term is positive:

$$\frac{P(1-\theta x)}{D}-1-\frac{\lambda(1-\theta)}{P_1} \geq 0 \tag{26}$$

We have:

$$\frac{1}{2}h[\frac{P(1-\theta x)}{D}-1-\frac{\lambda(1-\theta)}{P_1}] \times [(\frac{3-b_{p2}}{2-b_{p2}})(\frac{3-b_{p2}}{2-b_{p2}}-1)(\frac{\lambda(1-\theta)}{P_1}(1-b_{p2}))^{\frac{1}{2-b_{p2}}} tp_1^{\frac{3-b_{p2}}{2-b_{p2}}-2}$$
$$+(\frac{3-b_{p1}}{2-b_{p1}})(\frac{3-b_{p1}}{2-b_{p1}}-1)(1-b_{p1})^{\frac{1}{2-b_{p1}}} tp_1^{\frac{3-b_{p1}}{2-b_{p1}}-2}] \geq 0$$

$$\tag{27}$$

And the last positive word:

$$-\frac{1}{2}h[ts_1+ts_1(2)^{-b_s}] \times [(\frac{1}{2-b_{p2}})(\frac{1}{2-b_{p2}}-1)(\frac{\lambda(1-\theta)}{P_1}(1-b_{p2}))^{\frac{1}{2-b_{p2}}} tp_1^{\frac{1}{2-b_{p2}}-2}$$
$$+(\frac{1}{2-b_{p1}})(\frac{1}{2-b_{p1}}-1)(1-b_{p1})^{\frac{1}{2-b_{p1}}} tp_1^{\frac{1}{2-b_{p1}}-2}] \geq 0 \tag{28}$$

However the answer should be feasible, the Second derivation of the function should be positive, it means:

$$\frac{d^2TC(tp_1)}{dtp_1^2} = negetive\sec tion + \frac{1}{P_1-\lambda}(\frac{b_{p2}}{1-b_{p2}})(\frac{1-b_{p2}}{b_{p2}}-1)(\frac{\lambda(1-\theta)}{P_1})^{\frac{1-b_{p2}}{b_{p2}}} tp_1^{\frac{1-b_{p2}}{b_{p2}}-2} + positive\sec tion$$

Hence,

$$\frac{d^2TC(tp_1)}{dtp_1^2} = \{-h\frac{1}{(P-D-\lambda)}(\frac{1}{b_{p1}})tp_1^{-\frac{1}{b_{p1}}} + \frac{(1-b_{p1})^{\frac{1}{2-b_{p1}}}}{2-b_{p1}}(\frac{1}{2-b_{p1}}-1)ts_1(2)^{-b_s} tp_1^{\frac{1}{2-b_{p1}}-2} + \frac{h_1}{P_1}(-\frac{1}{1-b_{p1}})tp_1^{-\frac{1}{1-b_{p1}}-1}\} +$$

$$\{\frac{1}{P_1-\lambda}(\frac{b_{p2}}{1-b_{p2}})(\frac{1-b_{p2}}{b_{p2}}-1)(\frac{\lambda(1-\theta)}{P_1})^{\frac{1-b_{p2}}{b_{p2}}} tp_1^{\frac{1-b_{p2}}{b_{p2}}-2}\} + \{\frac{1}{\lambda}(\frac{1}{1-b_{p1}})tp_1^{-\frac{1}{1-b_{p1}}} + \frac{\lambda(1-\theta)}{P_1}(\frac{3-b_{p1}}{2-b_{p1}})(\frac{3-b_{p1}}{2-b_{p1}}-1)$$

$$(1-b_{p1})^{\frac{1}{2-b_{p1}}} tp_1^{\frac{3-b_{p1}}{2-b_{p1}}-2} + \frac{1}{2}[\frac{P(1-\theta x)}{D}-1-\frac{\lambda(1-\theta)}{P_1}] \times [(\frac{3-b_{p2}}{2-b_{p2}})(\frac{3-b_{p2}}{2-b_{p2}}-1)(\frac{\lambda(1-\theta)}{P_1}(1-b_{p2}))^{\frac{1}{2-b_{p2}}} tp_1^{\frac{3-b_{p2}}{2-b_{p2}}-2} +$$

$$+(\frac{3-b_{p1}}{2-b_{p1}})(\frac{3-b_{p1}}{2-b_{p1}}-1)(1-b_{p1})^{\frac{1}{2-b_{p1}}} tp_1^{\frac{3-b_{p1}}{2-b_{p1}}-2}] - \frac{1}{2}[ts_1+ts_1(2)^{-b_s}] \times [(\frac{1}{2-b_{p2}})(\frac{1}{2-b_{p2}}-1)(\frac{\lambda(1-\theta)}{P_1}(1-b_{p2}))^{\frac{1}{2-b_{p2}}}$$

$$tp_1^{\frac{1}{2-b_{p2}}-2} + (\frac{1}{2-b_{p1}})(\frac{1}{2-b_{p1}}-1)(1-b_{p1})^{\frac{1}{2-b_{p1}}} tp_1^{\frac{1}{2-b_{p1}}-2}]\} \geq 0$$

(29)

## 5. Numerical example

We consider 12000 units per year for production rate, 6000 units per year for demand rate, 2400 units per year for production rate of defective. Portion of produced items is assumed to be defective is 0.2, among the defective items, a portion $\theta=0.1$ is considered to be scrap and the other portion can be reworked and repaired.

As it seems:

$$\frac{P(1-\theta x)}{D}-1-\frac{\lambda(1-\theta)}{P_1} \geq 0 \Rightarrow \frac{12000(1-.2\times.1)}{6000}-1-\frac{2400(1-.1)}{8000} = .69 \geq 0$$

Other related parameters are summarized as follows:

M     $600 repair cost for breakdown,

K     $500 for each setup,

| | | |
|---|---|---|
| c | $4 per time, | |
| $c_r$ | $1 repaired cost for each item reworked, | |
| $c_s$ | $0.5 disposal cost for each scrap item, | |
| h | $1 per item per unit time, | |
| $h_1$ | $1.2 per item reworked per unit time, | |
| $ts_1$ | $0.02 years, time needed to repair and restore the machine. | |

Using the above data and based on the Fix learning rate, various figures are obtained by using Matlab software which shows the importance of learning effect in this type of EPQ model. Also by considering different Fix learning rate, (except type 1) cost function does not greatly differ.

In Figure 3 we suppose that all learning rates are the same and equal to 0.6 for both of perfect and defective products and set-up time. Minimum of this function is obtained 2.95.

$$\alpha_1, \alpha_2, \alpha_s = .6 \Rightarrow \frac{d^2 TC(tp_1)}{dtp_1^2} \geq 0$$

Please insert Figure 3 here.

In Figure 4 we suppose that all learning rates are the same and equal to 0.8. Minimum of this function is obtained 0.15 and this amount indicates that however learning rate increases the cost of production and time of production decreases which is better.

$$\alpha_1, \alpha_2, \alpha_s = .8 \Rightarrow \frac{d^2 TC(tp_1)}{dtp_1^2} \geq 0$$

Please insert Figure 4 here.

In Figure 5 we suppose that all learning rates are different in this condition, time of production increases (because of that all learning-rates are not increased). Minimum of this function is obtained 3.29.

$$\alpha_1 = .8, \alpha_2, \alpha_s = .6 \Rightarrow \frac{d^2 TC(tp_1)}{dtp_1^2} \geq 0$$

Please insert Figure 5 here.

We conclude from Figure 3, Figure 4 and Figure 5, learning rate have significant impact on the optimal point of production time and cost of production in this model of EPQ. If we increase the impact of learning (learning rate) optimal point of production time has been decreased.

**6. Conclusions and future research directions**

There are many key assumptions of the classical EPQ model which are inevitable in many manufacturing systems. Therefore, the study presents an inventory-production model with imperfect product, reworking, scrap, setup time, and machine breakdown to overcome the shortcomings of the traditional EPQ models. Moreover, learning effect in setup time and production time is considered. Based on a breakdown in machine, production-inventory cost function including breakdown cost, setup cost, non-defective products production cost, reworking cost and disposal cost of scrap items, and the holding cost of the non-defective and defective products is introduced. In the developed model, it is supposed that the machine has a fixed defective rate and a fixed scrap rate, and new mathematical model of production-inventory cost function has been introduced. Finally, a numerical example of the proposed model is presented which illustrates increasing the impact of learning (learning rate) leads to decrease optimal point of production time.

As a suggestion for future research we can extend this model by considering multi-product environment and use heuristic and meta-heuristic methods for solving the resulted NP-hard problem. In this paper, the effect of learning on setup and production times has been just discussed, so applying the forgotten effects could be considered as a future research.

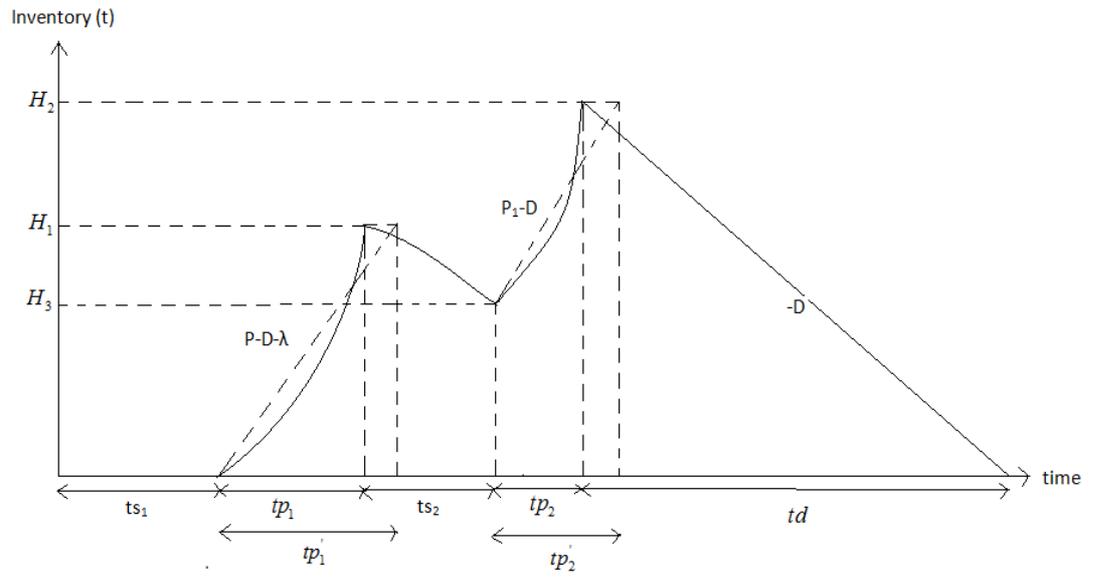

Figure 1. On-hand inventory of perfect quality in the developed EPQ model

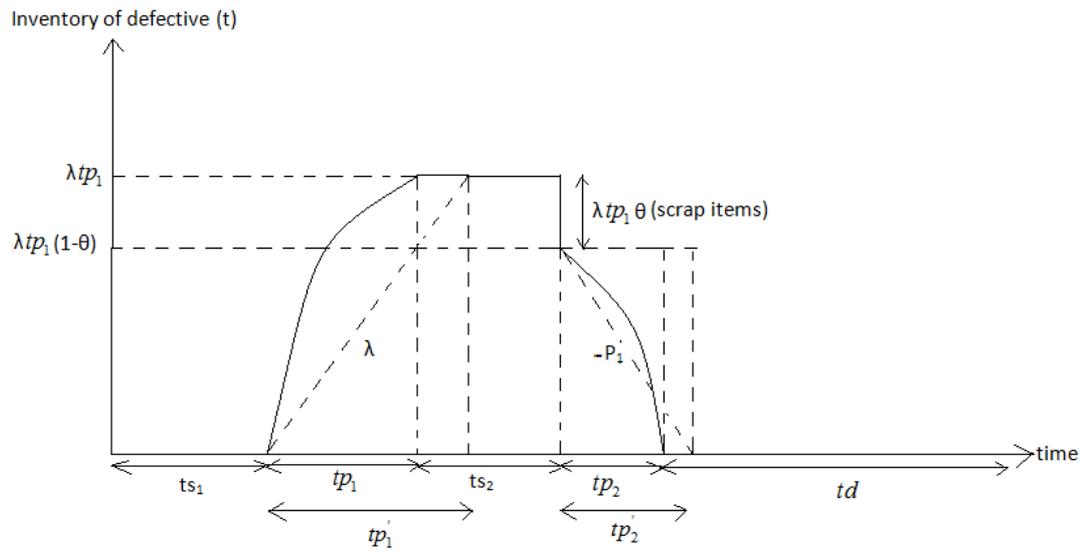

Figure 2. On-hand inventory of defected items (including scrap items) in the developed EPQ model

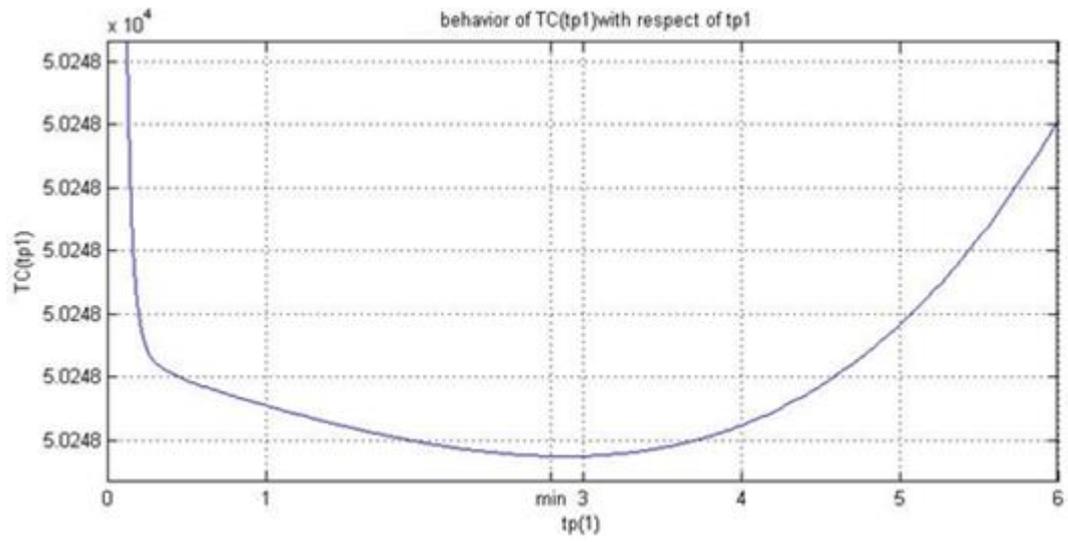

Figure 3. Behavior of TC(tp1) with respect to tp1

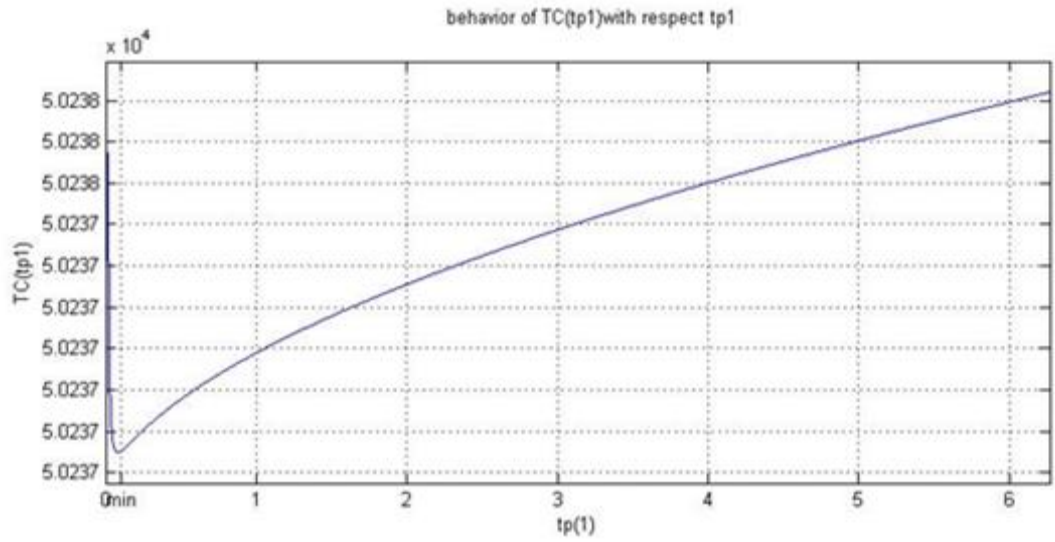

Figure 4. Behavior of TC(tp1) with respect to tp1

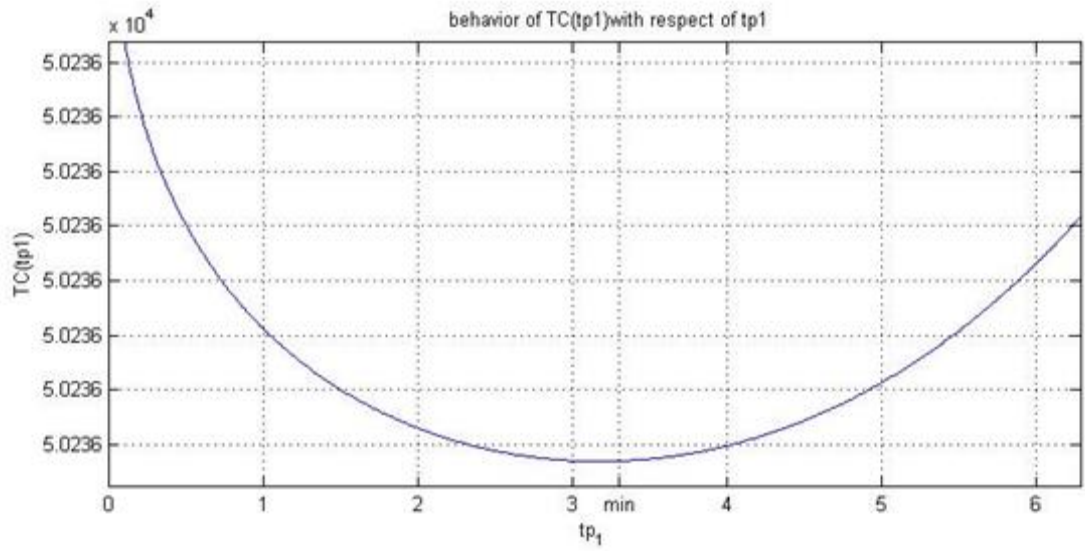

Figure 5. Behavior of TC(tp1) with respect to tp1